\documentclass[twoside,12pt]{article}
\usepackage{amsmath,graphicx}
\usepackage{amssymb}  
\newtheorem{theorem}{Theorem}[section]

\newtheorem{proposition}{Proposition}[section]

\newtheorem{remark}{Remark}[section]

\newtheorem{corollary}{Corollary}[section]

\def\Sum{\displaystyle\sum}

\def\sup{\displaystyle\sup}

\begin{document}

\begin{center}
{\bf   A free boundary problem for a discontinuous semilinear elliptic equations in convex ring .  }\vskip 0.5cm {\bf By} \vskip
1cm {\bf Sabri BENSID. }\\
\
\\ Dynamical Systems and Applications Laboratory\\ Department of Mathematics, Faculty of Sciences,\\ University Abou Bakr Bekaid of Tlemcen,\\ B.P.119, Tlemcen 13000,
Algeria.\\  Mail: $ edp\_sabri@yahoo.fr$

\end{center}
\begin{abstract}
In this paper, we consider the following  free boundary  problem $$
(P)\left\{\begin{array}{ll} \Delta u = \lambda \phi(x)\Sum_{i=1}^n  H(u-\mu_i )& \quad \mbox{
in }\ \Omega=\Omega_2\setminus \overline{\Omega}_1,
\\[0.3cm]u =0 &\quad \mbox{ on } \partial \Omega_2,
\\[0.3cm]u =M &\quad \mbox{ on } \partial \Omega_1.
\end{array}
\right. $$ The domain $\Omega$ is a convex ring where $\Omega_1$ and $\Omega_2$ are bounded convex domain in $\mathbb{R}^N,$ $N\geq 2,$
$H$ is the Heaviside step function, $\lambda, M,\mu_i$ $(i=1,..,n)$ are a given positive real parameters and $\phi$ is a given function . We show that under suitable conditions, there exists a solution to problem $(P)$ with a convex level set. We derive also an interesting result for the convexity of the set which is delimited by the free boundary in obstacle problem when the nonlinearity is discontinuous. At last, a detailed analysis is given by a construction of function $\phi$ which show that there are more solutions.
\end{abstract}

\section{Introduction and main results}
This paper is concerned with some geometrical properties of solutions of semilinear elliptic equation
$$\Delta u=f(u)\quad \hbox{in}\quad \Omega,$$
where $\Omega$ is a bounded domain of $\mathbb{R}^N$ and $f$ is a discontinuous nonlinearity. More precisely, let us be given a bounded domain $\Omega:=\Omega_2\setminus \overline{\Omega}_1$ where $\Omega_1,\Omega_2$ are bounded convex domains in $\mathbb{R}^N,$ $N\geq 2,$ ($\Omega_1\subset\subset\Omega_2$) and we study the following problem
$$
\left\{\begin{array}{ll} \Delta u = f(u)& \quad \mbox{
in }\ \Omega=\Omega_2\setminus \overline{\Omega}_1,
\\[0.3cm]u =0 &\quad \mbox{ on } \partial \Omega_2,
\\[0.3cm]u =M &\quad \mbox{ on } \partial \Omega_1,
\end{array}
\right. $$
where   $f$ is a discontinuous function with respect to $u$ and $M$ is a positive constant.\\
Our objective is to establish the convexity of the set $\{x\in \Omega,\quad u(x)\geq t\}$ for all $t\in (0,M)$ under a special nonlinearity function ( discontinuous).\\
Recall that when the function $f$ is continuous and nondecreasing, Caffarelli and Spruck in \cite{Caffspruck} prove the following result
\begin{theorem}\cite{Caffspruck}
Let $u$ be the unique solution of $$
\left\{\begin{array}{ll} \Delta u = \gamma(u)& \quad \mbox{
in }\ \Omega=\Omega_2\setminus \overline{\Omega}_1,
\\[0.3cm]u =0 &\quad \mbox{ on } \partial \Omega_2,
\\[0.3cm]u =M &\quad \mbox{ on } \partial \Omega_1,
\end{array}
\right. $$
with $\gamma(u)$ continuous and nondecreasing in $u,$ $\gamma(0)=0.$ Then, the level surfaces of $u$ are convex $C^{1+\alpha}$ hypersurfaces.
\end{theorem}
Later, many papers are devoted to study the convexity of superlevel set of the solutions of problem $(P)$ under weaker conditions. For instance, in \cite{Laurence1}, Laurence and Stredulinsky prove the existence of weak solutions with convex level lines in $\mathbb{R}^2$ of problem $(P)$ when $f$ satisfies $$f\in L^1(-\infty,+\infty),\quad f(x)\geq 0,\quad \hbox{and}\quad f(0)=0\quad \hbox{on}\quad (-\infty,0).$$
For many generalizations in the literature, we refer to \cite{Acker1},\cite{cafffrie}, \cite{Diaz1}, \cite{Friedmanphilips} and \cite{Laurence1} for further results.\\
In this work, we are concerned about study some qualitative properties of level set when $f$ is the sum of discontinuous nonlinearities i.e.
$$f(u)=\Sum_{i=1}^n  H(u-\mu_i ),\quad \mu_i\in (0,M),\quad i=\overline{1,n},$$
where $H$ is the Heaviside step function i.e
$$
H(t)=\left\{\begin{array}{ll} 1& \quad \mbox{
if }\ t\geq 0,
\\[0.3cm]0 &\quad \mbox{ if } t<0.
\end{array}
\right. $$
So, our first problem is
\begin{equation}
\label{first1}
\left\{\begin{array}{ll} \Delta u =\Sum_{i=1}^n  H(u-\mu_i ) & \quad \mbox{
in }\ \Omega=\Omega_2\setminus \overline{\Omega}_1,
\\[0.3cm]u =0 &\quad \mbox{ on } \partial \Omega_2,
\\[0.3cm]u =M &\quad \mbox{ on } \partial \Omega_1.
\end{array}
\right.
\end{equation}
The nonlinearity  $f$ is motivated by various problem of mathematical physics where many thresholds can appear. This phenomenon arise in several different context like fluid dynamics, potential flow in fluid mechanics, combustion theory. See for instance \cite{Acker1},\cite{Acker2},\cite{Bensidkaid}, \cite{cafffrie},\cite{Greco},\cite{henrot1},\cite{henrot2} and the references therein.\\\\
An interesting application is related to the growth tumor model. We expand this in more details.\\
Let us denote by $\Omega(t)\subset \mathbb{R}^N,$ $( N\geq 3),$ the spherical tumor region $ t>0$ and by $R(t)$ its radius. The classical diffusion
process shows that the oxygen concentration $u(x,t)$ satisfies the following equation, for $t>0$ and $x\in \Omega(t)$,$$
\frac{\partial u}{\partial t}=D \Delta u-\lambda f(u),
$$
with $D>0$ is the diffusion coefficient and $\lambda$ is a positive constant where $\lambda f(u)$ describe the consumed oxygen rate. Different
forms of the nonlinearity $f$ are given in the literature. See \cite{byrne}. A more general form of $f$ can be taken as
$$f(u)=1+\Sum_{i=1}^{n}\alpha_i H(\mu_i-u),$$
where   $\alpha_i\in \mathbb{R}$ for $i=1,2,...,n,$ $\mu_i$ is a critical value verifying $\mu_{i+1}>\mu_i>0.$ See \cite{bensid chekroun}. As in \cite{byrne}, we can consider $u$ quasi stable (the oxygen diffusion time scale
is much shorter than a typical tumor doubling time). So , our problem $(\ref{first1})$ provides a simple description of vascular tumor growth involves including a distributed source of nutrient. In this case, the developed tumor can be divided in several regions, an outer proliferating surrounds a quiescent annulus and a central necrotic core. For more details, see \cite{byrne},\cite{byrne1} and \cite{byrne2}.\\
In order to simplify the analysis, we fixed $t=t^*$ and thus $\Omega(t^*):=B_R\setminus B_{r_0},$ where $B_{\eta}$ is the ball of radius $\eta$ centred at $0.$ We can obtain an explicit information of model if we set $r=|x|$ where $u(r)$ verifies
\begin{equation}
\label{pb radial}
\left\{\begin{array}{ll} \Delta u(r) = \Sum_{i=1}^n  H(u(r)-\mu_i ) & \quad \mbox{
for }\ r_0<r<R,
\\[0.3cm]u(r_0) =M, &\quad \mbox{  }
\\[0.3cm]u(R) =0. &\quad \mbox{  }
\end{array}
\right.
\end{equation}
Now, for $s>0,$ let
$$\Omega_s^+(u):=\{x\in \Omega, u(x)\geq s\},\quad \Omega_s^-(u):=\{x\in \Omega, u(x)<s\},$$
$$\hbox{and}\quad \Gamma_s(u):=\{x\in \Omega, u(x)=s\}.$$
By a solution of problem $(1)$, we mean a function $u\in C^{1,\alpha}(\overline{\Omega}), 0<\alpha<1$  verifies $(1)$ such that the sets $\Gamma_{\mu_i}(u)$ are  analytic  hypersurfaces. The first main result of this work is the following theorem
\begin{theorem}
Let $M>0$  be a given positive constant and $\mu_1<\mu_2<...<\mu_n.$ Then  there exists a solution $u$ of problem $(\ref{first1})$ such that the set $\Omega_s^+(u)$ are convex for $s\in (0,M).$\\
In particular, the free boundaries $\Gamma_{\mu_i}(u)$ are  analytic convex hypersurfaces .
\end{theorem}
The result of Theorem $(1.2)$ gives rise a positive answer to the convexity of the unknown set $\Lambda$ when the nonlinearity is discontinuous in the following free boundary problem (obstacle problem):\\\\
Consider a solution $u$ of
\begin{equation}
\label{obstacle}
\left\{\begin{array}{ll} \Delta u =g(u) & \quad \mbox{
in }\ \Omega\setminus \Lambda,
\\[0.3cm]u =u_0 &\quad \mbox{ on } \partial \Omega,
\\[0.3cm]u =c &\quad \mbox{ on } \partial \Lambda,
\\[0.3cm]\partial_n u =0 &\quad \mbox{ on } \partial \Lambda,
\end{array}
\right.
\end{equation}
where $u_0$ is a given nonnegative constant, $\Lambda$ is a closed subset of $\Omega\subset \mathbb{R}^2$ (unknown region), $\partial_n u$ is the normal derivative of $u$ ( $n$ is the normal vector to $\partial \Lambda$) and $g$ is a given discontinuous function. In fact,  an immediate consequence is the following corollary.
\begin{corollary}
Assume that $g(u)=f_0+(1-f_0)H(u-\mu)$ where $f_0\in (0,1)$ and $u_0>c>0.$ If $\Omega$ is convex and if $u$ is
a solution of $(3)$, then $\Lambda$ is also convex.
\end{corollary}
Remark that the above corollary has been proved only in the case where $g$ is $C^0$ ( see \cite{Dolbeault}), so here, we give a weaker version.
The choice of function $g$ is motivated by the papers  \cite{B_D1} and $\cite{B_D2}.$\\\\
To the best of our knowledge, this is the first time in the literature that
problem $(3)$ is considered when the nonlinearity is discontinuous. There are two free boundaries to studied, first the free boundary  $\partial \Lambda$ (obstacle problem), second, the free boundary $\{x\in \Omega,\quad u(x)=\mu\}$  obtained by the form of discontinuous nonlinearity.\\\\
The last result of this work shows that we can find  an explicit construction to prove that multiple solutions can exists. We have the following theorem
\begin{theorem}
If $\Omega\subset \mathbb{R}^2,$ then for any $m\geq 2,$ there exists a $\phi\in L^{\infty}(\Omega)$ that admit at least $n m$ distinct solutions to
\begin{equation}
\label{first}
\left\{\begin{array}{ll} \Delta u =\lambda \phi(x)\Sum_{i=1}^n  H(u-\mu_i ) & \quad \mbox{
in }\ \Omega=B_R\setminus B_{r_0},
\\[0.3cm]u =0 &\quad \mbox{ on } \partial B_R,
\\[0.3cm]u =M &\quad \mbox{ on } \partial B_{r_0}.
\end{array}
\right.
\end{equation}
\end{theorem}
Our paper is organized as follows: In section 2, we prove theorem 1.2 concerning the convexity of level set and also derive the corollary 1.1. Section 3 is devoted to the proof of Theorem 1.3 using a detailed analysis of construction of $\phi $ and finally, some useful comments are given.
\section{Proof of Theorem 1.2}
First, we will approximate the Heaviside function $H(t)$ by increasing continuous functions $H_{\varepsilon,\mu}(t)$ so that
$$
H_{\varepsilon,\mu}(t)=\left\{\begin{array}{ll} 0& \quad \mbox{
if }\ t< \mu_1,
\\[0.3cm]\frac{t-\mu_1}{\varepsilon_1} &\quad \mbox{ if } \mu_1\leq t\leq \mu_1+\varepsilon_1,\\[0.3cm]
1& \quad \mbox{if } \mu_1+\varepsilon_1<t< \mu_2,
\\[0.3cm]\frac{t-\mu_2}{\varepsilon_2}+1 &\mbox{ if }\mu_2\leq t\leq  \mu_2+\varepsilon_2,
\\[0.3cm]. &\mbox{ }
\\[0.3cm]i-1 &\quad \mbox{ if } \mu_{i-1}+\varepsilon_{i-1}< t<  \mu_i,
\\[0.3cm]\frac{t-\mu_i}{\varepsilon_i}+i-1 &\quad \mbox{ if } \mu_i\leq t\leq \mu_i+\varepsilon_i,
\\[0.3cm]. &\quad \mbox{  }
\\[0.3cm]. &
\\[0.3cm]n &\quad \mbox{ if } t> \mu_n+\varepsilon_n,
\end{array}
\right. $$
where $\varepsilon=(\varepsilon_1,\varepsilon_2,..,\varepsilon_n)$ et $\mu=(\mu_1,\mu_2,..,\mu_n)$ satisfying
$$0<\varepsilon_i<\mu_{i+1}-\mu_{i},\quad \hbox{for}\quad i=1,..,n$$
We denote by   $u_{\varepsilon}$ the solution of of the following problem
$$
(P_\varepsilon)\left\{\begin{array}{ll} \Delta u = H_{\varepsilon,\mu}(u)& \quad \mbox{
in }\ \Omega=\Omega_2\setminus \overline{\Omega}_1,
\\[0.3cm]u =0 &\quad \mbox{ on } \partial \Omega_2,
\\[0.3cm]u =M &\quad \mbox{ on } \partial \Omega_1.
\end{array}
\right. $$
Since, the function $H_{\varepsilon,\mu}(t)$ is continuous and nondecreasing with $H_{\varepsilon,\mu}(0)=0,$ then according to the theorem 1.2 of \cite{Caffspruck}, the set $\Omega_t^+(u_{\varepsilon})$ are convex while $\Gamma_t(u_{\varepsilon})$ are smooth convex hypersurfaces.\\\\
Now, we want to prove that $u_{ \varepsilon}\rightarrow u$ as $\varepsilon\rightarrow0$ $(\varepsilon_i\rightarrow 0)$ in $C^{1,\alpha}(K), 0<\alpha<1$ for an arbitrary compact $K\subset \Omega.$ \\ The function $H_{\varepsilon,\mu}$ is uniformly bounded i.e $|H_{\varepsilon}|\leq n,$ then by the standard argument, $u_{\varepsilon}\in C^{1,\alpha}(K)$  and $u_{\varepsilon}$ is uniformly bounded in $C^{1,\alpha}(K).$ i.e
$$\|u_{\varepsilon}\|_{1,\alpha}\leq C<  \infty.$$
Hence, we can extract a  convergent subsequence also called $u_{\varepsilon}$ which converges uniformly to a function $v\in C^{1,\alpha}(\Omega).$\\
So, we can further show that if $u_{\varepsilon}\rightarrow u $ in the norm $W^{1,2}(\Omega),$ then by the uniqueness of the limit, we have $u=v.$\\\\
If we multiply our equation in $(P_{\varepsilon})$ by $(u-u_{\varepsilon})$ and integrate, we obtain
$$\int_{\Omega}(u-u_{\varepsilon})\Delta (u-u_{\varepsilon})dx=\int_{\Omega} \left(\Sum_{i=1}^nH(u-\mu_i)-H_{\varepsilon,\mu}(u_{\varepsilon})\right)(u-u_{\varepsilon})dx$$
$$=\int_{\Omega} \left(\Sum_{i=1}^nH(u-\mu_i)-H_{\varepsilon,\mu}(u)\right)(u-u_{\varepsilon})dx+\int_{\Omega} \left(H_{\varepsilon,\mu}(u)-H_{\varepsilon,\mu}(u_{\varepsilon})\right)(u-u_{\varepsilon})dx.$$
Hence,
$$-\int_{\Omega}|\nabla(u-u_{\varepsilon})|^2dx=\int_{\Omega} \left(\Sum_{i=1}^nH(u-\mu_i)-H_{\varepsilon,\mu}(u)\right)(u-u_{\varepsilon})dx$$$$+\int_{\Omega} \left(H_{\varepsilon,\mu}(u)-H_{\varepsilon,\mu}(u_{\varepsilon})\right)(u-u_{\varepsilon})dx.$$
Because the function $H_{\varepsilon,\mu}$ is monotone, then $\int_{\Omega} \left(H_{\varepsilon,\mu}(u)-H_{\varepsilon,\mu}(u_{\varepsilon})\right)(u-u_{\varepsilon})dx\geq 0$ and we have
$$\int_{\Omega}|\nabla(u-u_{\varepsilon})|^2dx\leq \int_{\Omega} \left(\Sum_{i=1}^nH(u-\mu_i)-H_{\varepsilon,\mu}(u)\right)(u-u_{\varepsilon})dx$$
$$\leq  \int_{\bigcup\limits_{i=1}^n\{0<u<\varepsilon_i\}}|u-u_{\varepsilon}|dx\leq \alpha m\left(\bigcup\limits_{i=1}^n\{0<u<\varepsilon_i\}\right),$$
where $m(.)$ is the Lebsegue measure and $\alpha>0.$\\\\
So, when $\varepsilon\rightarrow 0,$ ($\varepsilon_i\rightarrow 0$), then $\nabla u_{\varepsilon}\rightarrow \nabla u$ in $L^2(\Omega).$\\
Using the fact that $u-u_{\varepsilon}=0$ on $\partial \Omega,$ then in the norm $W^{1,2}(\Omega),$ $u_{\varepsilon}\rightarrow u$ and thus $u_{\varepsilon}\rightarrow u$ in $C^1(K),$ for a compact $K$ in $\Omega.$
\newline
\newline
Now, to prove the convexity of $\Omega_s^+(u),$ we have the following result
\begin{proposition}
The set $\Omega_s^+(u)=\bigcap\limits_{\varepsilon>0}  \Omega_s^+(u_{\varepsilon})$ is convex for $t\in (0,M).$
\end{proposition}
\textit{ \textbf{ Proof of Proposition 2.1}} Let $\varepsilon<\widetilde{\varepsilon}.$ We define
$$w_{\varepsilon,\widetilde{\varepsilon}}:=\{x\in \Omega, \quad u_{\varepsilon}(x)>u_{\widetilde{\varepsilon}}(x)\}.$$
If the set $w_{\varepsilon,\widetilde{\varepsilon}}\neq \varnothing,$ then
$$\Delta (u_{\varepsilon}-u_{\widetilde{\varepsilon}})=\Delta u_{\varepsilon}-\Delta u_{\widetilde{\varepsilon}}$$
$$=H_{\varepsilon}(u_{\varepsilon})-H_{\widetilde{\varepsilon}}(u_{\widetilde{\varepsilon}})\geq 0\quad \hbox{in}\quad w_{\varepsilon,\widetilde{\varepsilon}} $$
$$u_{\varepsilon}-u_{\widetilde{\varepsilon}}=0\quad \hbox{on}\quad \partial w_{\varepsilon,\widetilde{\varepsilon}} .$$
Then, by the maximum principle,
$$u_{\varepsilon}-u_{\widetilde{\varepsilon}}\leq 0\Rightarrow u_{\varepsilon}\leq u_{\widetilde{\varepsilon}}.$$
This is a contradiction. So, if $\varepsilon<\widetilde{\varepsilon},$ then $u_{\varepsilon}< u_{\widetilde{\varepsilon}}.$\\\\
Now, according to the theorem 1.2 of \cite{Caffspruck}, the set $\Omega_s^+(u_{\varepsilon})$ are convex and consequently that the intersection of convex set is convex, then $\bigcap\limits_{\varepsilon>0}  \Omega_s^+(u_{\varepsilon})$ is convex. By the monotonicity of $u_{\varepsilon}$ and the pointwise convergence to $u$, we can conclude that $\Omega_s^+(u):=\bigcap\limits_{\varepsilon>0}  \Omega_s^+(u_{\varepsilon})$ is convex.$\Box$\\Finally, we address the free boundaries regularity. We have the following
\begin{proposition}
For $\mu_i\in (0,M),$ $i=1,...,n,$ the free boundaries $\Gamma_{\mu_i}(u)$ are convex hypersurfaces of class $C^{1,\alpha},$ $0<\alpha<1.$
\end{proposition}
\textit{ \textbf{ Proof of Proposition 2.2}}\\\\
Since $\Omega_s^+(u)$ is convex, then at each point $x_0\in \Gamma_s(u),$ the condition of the sphere at $x_0$ is satisfied where $u$ is subharmonic.\\ By the strong maximum principle for subharmonic functions, we have $|\nabla u(x_0)|\neq 0,$ $\forall x_0\in \Gamma_s(u).$ Hence, because $u\in C^{1,\alpha}(\Omega),$ $\alpha\in (0,1),$ then according to the implicit function theorem, we conclude that for $i=1,..,n$ the free boundaries $\Gamma_{\mu_i}$ are  hypersurface of class $C^{1,\alpha}.$ $\Box$\\\\
To conclude the proof of Theorem 1.2, it suffice to see that by the application of hodograph transformation method, we have the analyticity of $\Gamma_{\mu_i}.$ This method was used by the author in \cite{Bensid2}.  We refer the reader to \cite{Stamppachia}  for a complete description of this method.\\\\
Now, we will prove the corollary 1.1 in a similar way. Recall that Dolbeault and Monneau prove in \cite{Dolbeault} that if $\Omega\subset \mathbb{R}^2$ is convex and $u$ is solution of the following problem
\begin{equation}
\label{obstacle}
\left\{\begin{array}{ll} \Delta u =f(u) & \quad \mbox{
in }\ \Omega\setminus \Lambda,
\\[0.3cm]u =u_0 &\quad \mbox{ on } \partial \Omega,
\\[0.3cm]u =0 &\quad \mbox{ on } \partial \Lambda,
\\[0.3cm]\partial_n u =0 &\quad \mbox{ on } \partial \Lambda,
\end{array}
\right.
\end{equation}
where $f(0)>0$ and $u\longmapsto f(u)$ is increasing function of class $C^0,$ then $\Lambda$ is also convex. (See Theorem 1 in \cite{Dolbeault}).\\
Hence, we approximate the solution of problem $(3)$ by functions $u_{\varepsilon}$ such that
\begin{equation}
\label{obstacle}
\left\{\begin{array}{ll} \Delta u_{\varepsilon} =g_{\varepsilon}(u_{\varepsilon}) & \quad \mbox{
in }\ \Omega\setminus \Lambda_{\varepsilon},
\\[0.3cm]u_{\varepsilon} =u_0 &\quad \mbox{ on } \partial \Omega,
\\[0.3cm]u_{\varepsilon} =c &\quad \mbox{ on } \partial \Lambda_{\varepsilon},
\\[0.3cm]\partial_n u_{\varepsilon} =0 &\quad \mbox{ on } \partial \Lambda_{\varepsilon},
\end{array}
\right.
\end{equation}
where $$g_{\varepsilon}(t):=\left\{\begin{array}{ll} f_0 & \quad \mbox{
 }\ t<\mu,
\\[0.3cm](1-f_0)\frac{(t-\mu)}{\varepsilon} &\quad \mbox{  } \mu\leq t\leq \mu+\varepsilon,
\\[0.3cm]1 &\quad \mbox{  } t>\mu+\varepsilon.
\end{array}
\right.$$
Since $g_{\varepsilon}(u)$ is continuous, nondecreasing and verifies $g_{\varepsilon}(0)=f_0>0,$ then from theorem 1 in \cite{Dolbeault}, we assure that $\Lambda_{\varepsilon}:=\{u_{\varepsilon}>c\}$ is convex.\\
Then, by our result obtained in Theorem 1.2, we have $u_{\varepsilon}$ converge to the solution of problem $(3)$ in $C_{loc}^1$ as $\varepsilon\rightarrow 0.$ So, by proposition 2.1, the set $\Lambda:=\bigcap\limits_{\varepsilon>0} \Lambda_{\varepsilon}$ is convex.
\begin{remark}
In \cite{Dolbeault}, the authors consider the following general free boundary problem
\begin{equation}
\label{obstacle}
\left\{\begin{array}{ll} div(a (|\nabla u|^2) \nabla u) =f(u) & \quad \mbox{
in }\ \Omega\setminus \Lambda,
\\[0.3cm]u =u_0 &\quad \mbox{ on } \partial \Omega,
\\[0.3cm]u =0 &\quad \mbox{ on } \partial \Lambda,
\\[0.3cm]\partial_n u =0 &\quad \mbox{ on } \partial \Lambda,
\end{array}
\right.
\end{equation}
where $a$ is an increasing function of class $C^1.$ It is easy to adapt the result of Theorem 1 in \cite{Dolbeault} for problem $(6).$ ($a\equiv 1$ )
\end{remark}
\section{Proof of Theorem 1.3}
In this section, we will prove the theorem 1.3. We denote by $\Omega_{\rho,\eta}(X_0):=B_{\rho}\setminus B_{\eta}$ where $B_{\rho}$ and $B_{\eta}$ are ball of radius $\rho$ and $\eta$ centred at $X_0$ respectively.\\ Let $u$ be the function solving $$\left\{\begin{array}{ll} \Delta u =\lambda \Sum_{i=1}^n  H(u-\mu_i ) & \quad \mbox{
in }\ \Omega_{R,r_0},
\\[0.3cm]u =0 &\quad \mbox{ on } \partial B_R,
\\[0.3cm]u =M &\quad \mbox{ on } \partial B_{r_0},
\end{array}
\right.$$
such that $\{u>\mu_i\}:=\Omega_{r_i,r_{i-1}}(0),$ for $i=1,2,...,n.$\\\\
Let $g$ be a conformal mapping such that our region $\Omega_{r_i,r_{i-1}}(0)$ is mapped onto $\Omega_{\delta_i,\delta_{i-1}}(X_0)$ where $X_0=(x_0,0)$ for some $x_0$ and $\delta_i$ . So, for a fixed $R,r_0$ the function $g$ maps $\Omega_{R,r_0}$ onto itself where $g^{-1}=g $ and $|g(z)|^2\leq M<\infty$ for $z\in \Omega_{R,r_0}(0)$ where $z=(x,y).$\\
We define a new function $v_1(z)$ by $$v_1(z)=u(g^{-1}(z))$$
verifying the equation $\Delta v_1=\Delta u|g'(z)|^2.$\\ Hence, $v_1$ satisfies
 $$\left\{\begin{array}{ll} \Delta v_1 =\lambda\phi_1 \Sum_{i=1}^n  H(v_1-\mu_i ) & \quad \mbox{
in }\ \Omega_{R,r_0},
\\[0.3cm]v_1 =0 &\quad \mbox{ on } \partial B_R,
\\[0.3cm]v_1 =M &\quad \mbox{ on } \partial B_{r_0},
\end{array}
\right.$$
where $\phi_1=|g'(z)|^2$ and $\{v_1>\mu_i\}=\{(x,y)\in \mathbb{R}^2, \delta_{i-1}\leq (x-x_0)^2+y^2\leq \delta_i\}:=S_{1,i}.$\\\\
Let $m\geq 2$ be given. Then, for $\delta_i,\delta_{i-1}$ small enough and under a rotation of angle $\frac{2\pi}{m},$ the image of $S_{1,i}$ does not intersect the original set $S_{1,i}.$ This images are denoted by $S_{k,i} $ under $k$ rotations. Hence, each rotation gives a new functions $v_k$ verifying $$\left\{\begin{array}{ll} \Delta v_k =\lambda\phi_k \Sum_{i=1}^n  H(v_k-\mu_i ) & \quad \mbox{
in }\ \Omega_{R,r_0},
\\[0.3cm]v_k =0 &\quad \mbox{ on } \partial B_R,
\\[0.3cm]v_k =M &\quad \mbox{ on } \partial B_{r_0},
\end{array}
\right.$$
where $\{v_k>\mu_i\}=S_{k,i}.$\\\\
Now, define a continuous function $\widetilde{\phi}$ by $$\widetilde{\phi}=\phi_k\quad \hbox{in}\quad \bigcup\limits_{i=1}^n S_{k,i}  $$
and we extend $\widetilde{\phi}$ in a continuous way to the rest of the domain $\bigcup\limits_{i=1}^n S_{k,i}.$ So, $v_k$ satisfies
$$\left\{\begin{array}{ll} \Delta v_k =\lambda\widetilde{\phi} & \quad \mbox{
in }\ S_{k,1},
\\[0.3cm]\Delta v_k =2\lambda\widetilde{\phi} &\quad \mbox{ in } S_{k,2},\\
.\\
.\\
\\[0.3cm]\Delta v_k =0 &\quad \mbox{ in } \Omega_{R,r_0}\setminus \bigcup\limits_{i=1}^n S_{k,i}.
\end{array}
\right.$$
With this $\widetilde{\phi},$ the function $v_k$ verifies
$$\left\{\begin{array}{ll} \Delta v_k =\lambda\widetilde{\phi} \Sum_{i=1}^n  H(v_k-\mu_i ) & \quad \mbox{
in }\ \Omega_{R,r_0},
\\[0.3cm]v_k =0 &\quad \mbox{ on } \partial B_R,
\\[0.3cm]v_k =M &\quad \mbox{ on } \partial B_{r_0}.
\end{array}
\right.$$
Then, we can conclude that problem $(\ref{first})$ has at least $n.m$ solutions.
\section*{Final remarks}
In this paper, we prove the convexity of the level set corresponding to semilinear elliptic equation involving discontinuous nonlinearities. It is interesting to generalize the same result ( Theorem 1.2) to elliptic equations in divergence form.\\\\
Also, it is a fruitful idea to study the asymptotic behavior of domains\\ $\{x\in \Omega, u(x)>\mu_i\}$ as $\lambda\rightarrow +\infty.$ This question can be very important in growth tumor models. \\\\
Finally, a question of regularity can be asked. When will the
free boundary of problem $(P)$ develop singularities? The characterization of free boundary
and his regularity remains an open problem precisely when the domain $\Omega$ is
not smooth enough.

\end{document}